\documentclass[a4paper,11pt]{article}
\usepackage{indentfirst,latexsym,bm,amsmath,amssymb,amsthm}
\usepackage[dvips]{graphicx}
\usepackage{color}

\makeatletter\@addtoreset{equation}{section} \makeatother
\newtheorem{theorem}{Theorem}[section]
\newtheorem{Lemma}{Lemma}[section]

\newtheorem{remark}{Remark}[section]

 

\title{A class of global large, smooth solutions for the magnetohydrodynamics with the Hall and ion-slip effects}

\begin{document}
	\maketitle
	
	\centerline{\scshape Huali Zhang$^*$}
	\medskip
	{\footnotesize
		\centerline{Changsha University
			of Science and Technology}
		\centerline{School of Mathematics and Statistics}
		\centerline{ Changsha, 410114, People's Republic of China.}
	} 
	
	\medskip
	
	
	\bigskip
	

	\begin{abstract}
		In this paper, the Cauchy's problem for fractional MHD system with the Hall and ion-slip effects is considered. By exploring the structure of semilinear and quasilinear terms, we prove the global existence of solutions for a class of large initial data. Both the velocity and magnetic fields could be arbitrarily large in $H^3(\mathbb{R}^3)$.
	\end{abstract}
	
	\section{Introduction}
	
	In this paper we consider the following incompressible magnetohydrodynamics with the Hall and ion-slip effects:	
	\begin{equation}\label{HMHD}
		\begin{cases}
			u_t+\nu \Lambda^{\alpha} u+u \cdot \nabla u+\nabla p - b \cdot \nabla b=0,\\
			b_t+ \mu \Lambda^{\beta} b+u \cdot \nabla b-b \cdot \nabla u+ \sigma \nabla \times \left( \left(\nabla \times b \right) \times b\right)- \kappa \nabla \times \left( \left((\nabla \times b) \times b \right) \times b\right)=0,\\
			\nabla\cdot u=0, \quad \nabla\cdot b=0,\\
			u|_{t=0}=u_0, \ b|_{t=0}=b_0,
		\end{cases}
	\end{equation}
	on the domain $(t,x) \in \mathbb{ R}^+ \times \mathbb{R}^3$, where $\Lambda=\sqrt{-\Delta}$, $\beta=2, \alpha \in [0,2]$. Here $u=(u_1,u_2,u_3)^{\text{T}}, b=(b_1,b_2,b_3)^{\text{T}}\in \mathbb{R}^3$ denote the fluid velocity and magnetic fields respectively. The scalars $p, \nu, \mu$ are the pressure, viscosity, magnetic diffusivity respectively ($\nu$, $\mu$ are positive constants). $\kappa \geq 0, \sigma$ are constants. The Hall term $ \nabla \times \left( \left(\nabla \times b \right) \times b\right)$ is for the Hall effect, and $ \nabla \times \left( \left(\nabla \times b \times b \right) \times b\right)$ for ion-slip effect. $u_0$ and $b_0$ are the initial data satisfying
	\begin{equation}\label{HMHD1}
		\nabla \cdot u_0=\nabla \cdot h_0=0.
	\end{equation}
	
	Equation \ref{HMHD} is important to describe some physical phenomena, e.g., in the magnetic reconnection in space plasmas, star formation, neutron stars and dynamo.
	In the case $\sigma=\kappa=0$,  Equation \ref{HMHD} reduces to the standard MHD equations; in the case $\kappa=0$, Equation \ref{HMHD} reduces to Hall-MHD system. They has been extensively researched by a lot of excellent works, for instances \cite{AD}--\cite{CL}, \cite{CWW}--\cite{D}, \cite{LLZ}--\cite{LZZ} and \cite{JO,ST,RWX,ZZ,YM}.	

		For MHD system with the Hall and ion-slip effects, there are some interesting results related to the well-posedness theory, see references \cite{M,MF,MS}. Recently, Fan et.al in \cite{FZ} established global existence and time decay for small solutions. Very recently, Zhao and Zhu in \cite{Zz} gave a proof of global existence for small solutions under weaker smallness conditions.
	However, none of results are known for MHD system with the Hall and ion-slip effects for general initial data without smallness conditions.  It's quite rare to prove the existence of large, smooth, global solutions for quasilnear system. Under a class of large initial data, we found some results for incompressible Navier-Stokes equations and incompressible standard MHD equations, see \cite{CM, CG, LLZ, LZZ, ZZ, LW, Z} for details. Those motivate us to study the global well-posedness of Cauchy's problem of Equation \ref{HMHD} with large inital data. But the Hall and ion-slip term heightens the level of nonlinearity of the standard MHD system from a second-order semilinear to a second-order quasilinear level, significantly making its qualitative analysis more difficult. To the author's knowledge, it's quite rare to prove the existence of large, smooth, global solutions for quasilnear system. Based on the idea of Lei-Lin-Zhou in \cite{LLZ}, we study the existence of global large solutions. It's not a trivial extension from Hall-MHD equations to Equation \ref{HMHD}, for the additional difficulty of Equation \ref{HMHD} arise from: (1) the cubic quasi-linear term which requires some new cancellation estimates; (2) the proof of nonlinear smallness for some cross term. Using the large initial data constructed in \cite{LLZ} and combining nonlinear structures with commutator energy estimates, we are very fortunately to go through these difficulties---estimates \eqref{f0}-\eqref{f1} and \eqref{21}.

	The aim of this paper is to prove the existence of a unique, global smooth solution of MHD system with the Hall and ion-slip effects in $H^3(\mathbb{R}^3)$. Our result completely drops the smallness condition on the initial data.

	Before we state our main results, we first give some notations. Let  $\chi(x) \in C^{\infty}_0(\mathbb{R}^3)$ be a cut off function satisfying $|\chi(x)| \leq 2$ and
	\begin{equation}\label{102}
		\begin{split}
			&\chi(x)\equiv1, \quad \text{for} \ |x| \leq 1; \quad \chi(x)\equiv 0, \quad \text{for} \ |x| \geq 2,
			\\
			& |\nabla^k \chi(x)| \leq 2, \quad 0 \leq k \leq 5.
		\end{split}
	\end{equation}
	Denote
	\begin{equation}\label{M0}
		\chi_{M_0}(x):=\chi(\frac{x}{M_0}).
	\end{equation}
	Here $M_0$ is a positive constant. Let $v_0$ be that constructed by Lei et al. \cite{LLZ}, and it has the following properties
	\begin{align}\label{103}
		&\nabla \cdot v_0=0, \quad \nabla \times v_0=\sqrt{-\Delta}v_0,
		\\
		& \text{supp} \hat{v}_0 \subseteq \{ \xi| 1-\delta \leq |\xi| \leq 1+\delta \}, \quad 0<\delta\ \leq \frac{1}{2},
		\\
		&||\hat{v}_0||_{L^1} \leq M_1, \quad |\nabla^k v_0| \leq \frac{M_2}{1+|x|}, \quad 0\leq k \leq 5,
	\end{align}
	where $M_1, M_2$ are positive constants, $\hat{v}_0$ is the Fourier transform of $v_0$ and the operator $\sqrt{-\Delta}$ is defined through the Fourier transform
	\begin{equation*}
		\widehat{\sqrt{-\Delta}f}(\xi)=|\xi|\hat{f}(\xi).
	\end{equation*}
	
	Our main result is as follows.
	\begin{theorem}\label{thm}
		Consider Cauchy's problem \ref{HMHD}-\ref{HMHD1}. Suppose that  	
		\begin{align}\label{104}
			&\quad u_0=u_{01}+\chi_{M_0} u_{02},
			\\
			&\quad b_0=b_{01}+\chi_{M_0} b_{02}.
		\end{align}	
		with
		\begin{align}\label{100}
			&\nabla \cdot u_{02}=\nabla \cdot b_{02}=0,
			\\
			&u_{02}=\alpha_1 v_0, \quad b_{02}=\alpha_2 v_0,
		\end{align}
		where $\chi_{M_0},v_0$ are stated as above. $\alpha_1,\alpha_2$ are two real constants. Then there exist constants $\delta^{-\frac{1}{2}} \geq M_0 \gg 1$ depending on $M_1, M_2, \alpha_1, $ $\alpha_2, \mu, \nu, \eta$ such that Cauchy's problem \ref{HMHD}-\ref{HMHD1} has a unique, global smooth solution provided that
		\begin{equation}\label{111}
			||u_{01}||_{H^3}+||b_{01}||_{H^3} \leq M_0^{-\frac{1}{2}}.
		\end{equation}
	\end{theorem}
	\begin{remark}
		For
		\begin{equation*}
			\begin{split}
				&	||u_0||_{L^\infty}+||b_0||_{L^\infty} \leq  M_1,
				\\
				&	||u_0||_{H^3}+||b_0||_{H^3}\leq  \left( M_0^{-\frac{1}{2}}+(|\alpha_1|+|\alpha_2|)\sum_{k=0}^{3}\frac{M_2}{M_0^k} \right),
			\end{split}
		\end{equation*}
		and the constant $M_1,M_2$ can be arbitrary large, thus our initial data can be arbitrary large. Comparing with \cite{FZ,Zz}, our result can be seen as an non-trival improvement of Fan et al.'s and Zhao et al.'s work, for we completely drops the smallness condition on the initial data. 
	\end{remark}
\begin{remark}
	The parameter $\alpha, \beta$ indicates the strength of dissipation for velocity and magnetic field respectively. If the parameter $\alpha, \beta$ is larger, then the corresponding dissipation is stronger. When $\sigma=\kappa=0$, the conclusion in Theorem 1.1 still holds for all $\alpha, \beta \in [0,2]$. If $\sigma \neq 0$ or $\kappa \neq 0$, considering the quasilinear terms for magnetic field on \eqref{HMHD},  then the strong dissipative term $\mu \Delta b  (\beta=2)$ may be necessary to compensate for the loss of regularity in exploring large solutions.
\end{remark}
	\begin{remark}
		In the limiting case $\delta =
		0$, $\nabla \times u_{02} = u_{02}$, $\nabla \times b_{02} = b_{02}$, and the flow, magnetic field are called Beltrami flow and force-free fields respectively. Let us also mention that the magnetic energy achieves the minimum value for force-free fields, one can refer \cite{T} for details.
	\end{remark}
	\begin{remark}
		We throughout use a notation $C$. It may be different from line to line, but it is a universal positive constant in this paper.
	\end{remark}
	The proof of Theorem \ref{thm} is based on a perturbation argument along with a standard cut-off technique, and the perturbation is as large as the initial data. Compared with Hall-MHD equations, a part of the nonlinearities may not be small for Equation \ref{HMHD} (see \eqref{fg}). Fortunately, by combining the nonlinear structure of the term and commutator estimates, these terms can be estimated carefully.

	This paper is organized as follows: In section 2, we introduce commutator estimates and give some estimate of quadratic and cubic terms. Section 3 is devoted to prove the global existence and uniqueness of large smooth solutions for Equation \ref{HMHD}.
	
	\section{Preliminaries}
	\begin{Lemma}\label{lemc}\cite{CV}
		Let $s>0$. Let $p, p_2, p_3 \in(1,\infty)$ and $p_2,p_4 \in[1,\infty)$ satisfy
		\begin{equation*}
		\frac{1}{p}=\frac{1}{p_1}+\frac{1}{p_2}=\frac{1}{p_3}+\frac{1}{p_4}.
		\end{equation*}
		Then there exist two constants $C_1,C_2$,
		\begin{equation*}
		\begin{split}
		&|| \Lambda^s(fg)||_{L^p} \leq C_1 \left( || \Lambda^s f||_{L^{p_1}} ||g||_{L^{p_2}}+ || \Lambda^{s}g||_{L^{p_3}} ||f||_{L^{p_4}} \right),
		\\
		&|| [\Lambda^s,f] g||_{L^p} \leq C_2 \left( || \Lambda^s f||_{L^{p_1}} ||g||_{L^{p_2}}+ || \Lambda^{s-1}g||_{L^{p_3}} ||\nabla f||_{L^{p_4}} \right).
			\end{split}
		\end{equation*}
	\end{Lemma}
	Let $f,g$ satisfy
	\begin{equation}\label{300}
		\begin{cases}
			f_t+\nu \Lambda^\alpha f=0,\\
			t=0: f=u_{02},
		\end{cases}
	\end{equation}
	and
	\begin{equation}\label{301}
		\begin{cases}
			g_t-\mu \Delta g=0,\\
			t=0: g=b_{02}.
		\end{cases}
	\end{equation}
	Therefore, we have
	\begin{equation*}
		f=e^{-\nu t \Lambda^\alpha}u_{02}, \quad g=e^{\mu t\Delta}b_{02}.
	\end{equation*}
	\begin{Lemma}\label{tuilun1}
		Let $f,g$ be defined in \eqref{300}, \eqref{301}. It holds
		\begin{align*}
			&\nabla \cdot f=0, \quad \quad \nabla \times f=\sqrt{-\Delta} f,
			\\
			& \nabla \cdot g=0, \quad \quad \nabla \times g=\sqrt{-\Delta} g,
			\\
			&|\nabla^k f| \leq \frac{|\alpha_1| M_2}{1+|x|}e^{-\frac{\nu t}{2^\alpha}}, \quad |\nabla^k g| \leq \frac{|\alpha_2| M_2}{1+|x|}e^{-\frac{\mu t}{4}}, \quad 0\leq |k| \leq 5.
		\end{align*}
	\end{Lemma}
\begin{proof}
	By
	\begin{align*}
	& \nabla \cdot v_0=0, \quad \nabla \times v_0=\sqrt{-\Delta}v_0,
	\\
	& f=e^{-\nu t\Lambda^\alpha}u_{02}, \quad g=e^{\mu t\Delta}b_{02},
	\end{align*}
	we can deduce that
	\begin{align*}
	&\nabla \cdot f=0, \quad \quad \nabla \times f=\sqrt{-\Delta} f,
	\\
	& \nabla \cdot g=0, \quad \quad \nabla \times g=\sqrt{-\Delta} g.
	\end{align*}
	We choose a $C^\infty(\mathbb{R}^3)$ cut-off function $\gamma(\xi)$ such that $\alpha\equiv1$ on the support of $v_0$, and $\gamma(\xi)\equiv0$ if $|\xi|\geq 1+2\delta$ or $|\xi|\leq 1-2\delta$. Then we have
	$$
	f(t,x)=\alpha_1e^{-\frac{\nu t}{2^\alpha}} {\mathcal{F}}^{-1}\left(e^{-\nu(|\xi|^\alpha-\frac{1}{2^\alpha})t}\gamma(\xi)\right) \ast v_0,
	$$
	$$
	g(t,x)=\alpha_2e^{-\frac{\mu t}{2}} {\mathcal{F}}^{-1}\left(e^{-\mu(|\xi|^2-\frac{1}{2})t}\gamma(\xi)\right) \ast v_0.
	$$
	In a result, we get
	$$|\nabla^k f| \leq \frac{|\alpha_1| M_2}{1+|x|}e^{-\frac{\nu t}{2^\alpha}}, \quad |\nabla^k g| \leq \frac{|\alpha_2| M_2}{1+|x|}e^{-\frac{\mu t}{4}}, \quad 0\leq |k| \leq 5.$$
\end{proof}
	\begin{Lemma}\label{tuilun2}
		Set $\tilde{f}:=\chi_{M_0}f, \tilde{g}:=\chi_{M_0}g$. Let $f,g, \chi_{M_0}$ be defined by \eqref{300}, \eqref{301} and \eqref{M0} respectively. Then we have
		\begin{equation}\label{0000}
			||\tilde{f}||_{W^{5,\infty}}+||\tilde{g}||_{W^{5,\infty}} \leq C\left(|\alpha_1|M_1 e^{-\frac{\nu t}{2^\alpha}}+|\alpha_2|M_1 e^{-\frac{\mu t}{4}}\right),
		\end{equation}
		\begin{equation}\label{0001}
			\begin{split}
				& ||\tilde{f} \times \left( \nabla \times \tilde{f}\right)||_{H^3}
				+||\tilde{g} \times \left( \nabla \times \tilde{g}\right)||_{H^3}
				\\
				\leq
				& C\left( \alpha_1^2 e^{-\frac{\nu t}{2^{\alpha-1}}}+\alpha_2^2 e^{-\frac{\mu t}{2}}\right)\left( \delta M_0^{\frac{3}{2}}M_1^2+M_0^{-1}M_2^2\right),
			\end{split}
		\end{equation}
		\begin{equation}\label{0003}
		||\left( (\nabla \times \tilde{g}) \times \tilde{g} \right) \times \tilde{g} ||_{H^3} \leq C\left( \delta M_0^{\frac{3}{2}}M_1^3+M_0^{-1}M_2^3 \right)|\alpha_2|^3 e^{-\frac{3\mu t}{4}},
		\end{equation}
		\begin{equation}\label{0002}
			\int^\infty_0 || \tilde{f} \times \tilde{g}||_{H^3}(t)dt \leq CM_0^{\frac{3}{2}}M_1^2 (1+\alpha)\delta.
		\end{equation}
	\end{Lemma}
	\begin{proof}
		Firstly, we have $\left|\nabla^k \chi_{M_0} \right| \leq  CM^{-k}_0, k\leq 5$. Then
		\begin{align}\label{200}
		||\tilde{f}||_{W^{5,\infty}}&=||\chi_{M_0}f||_{W^{5,\infty}} \leq ||\chi_{M_0}||_{W^{5,\infty}}||f||_{W^{5,\infty}}
		\leq C ||f||_{W^{5,\infty}}.
		\end{align}
		Using $\hat{f}=e^{-\nu  |\xi|^\alpha t} \hat{u}_{02}$ and $ \text{supp}  \ \hat{u}_{02} \subseteq \left\{ \xi| 1-\delta \leq |\xi| \leq 1+\delta \right\}$, $0 < \delta \leq \frac{1}{2}$, we get
		\begin{align*}
		||f||_{W^{5,\infty}} \leq  \left| \left| (1+|\xi|)^5 \hat {f}   \right| \right|_{L^1_\xi}  \leq C \left| \left|  e^{-\nu t |\xi|^\alpha} \hat{u}_{02}  \right| \right|_{L^1_\xi}  \leq C |\alpha_1| M_1 e^{-\frac{\nu t} {2^\alpha}}.
		\end{align*}
		Similarly, we have
		\begin{align}\label{201}
		||g||_{W^{5,\infty}} \leq C |\alpha_2| M_1 e^{-\frac{\mu t} {4}}.
		\end{align}
		Adding \eqref{201} to \eqref{200}, we obtain
		\begin{align*}
		||\tilde{f}||_{W^{5,\infty}}+||\tilde{g}||_{W^{5,\infty}} \leq C\left(|\alpha_1| M_1 e^{-\frac{\nu t}{2^\alpha}}+|\alpha_2|M_1 e^{-\frac{\mu t}{4}}\right).
		\end{align*}
		Secondly, we notice the fact
		$$\nabla \times (\chi_{M_0}f)=\nabla \chi_{M_0} \times f +\chi_{M_0} \nabla \times f,$$
		$$\nabla \times (\chi_{M_0}g)=\nabla \chi_{M_0} \times g +\chi_{M_0} \nabla \times g.$$
		Thus, we get
		\begin{align*}
		&|| \tilde{f} \times (\nabla \times \tilde{f})||_{H^3} + || \tilde{g} \times (\nabla \times \tilde{g})||_{H^3}
		\\
		=&|| \chi_{M_0}f \times \left(\nabla \times (\chi_{M_0}f)\right)||_{H^3}+|| \chi_{M_0}g \times \left(\nabla \times (\chi_{M_0}g)\right)||_{H^3}
		\\
		\leq &C || \chi^2_{M_0}||_{H^3} \left(  ||f \times (\nabla \times f)||_{W^{3,\infty}} + ||g \times (\nabla \times g)||_{W^{3,\infty}}\right)
		\\
		& \quad \quad + C||\nabla(\chi^2_{M_0})||_{W^{3,\infty}} \left(  || |f|^2 ||_{H^3} + || |g|^2||_{H^3}\right)
		\end{align*}
		and
		\begin{equation}\label{g2}
		\begin{split}
		&||\left( (\nabla \times \tilde{g}) \times \tilde{g} \right) \times \tilde{g} ||_{H^3}
		\\
		= & ||\left( (\nabla \times (\chi_{M_0}g)) \times (\chi_{M_0}g) \right) \times (\chi_{M_0}g) ||_{H^3}
		\\
		 \leq & C \left(|| \chi^3_{M_0}||_{H^3} ||\left( (\nabla \times g) \times g \right) \times g ||_{W^{3,\infty}} + ||\nabla \chi_{M_0}^3||_{W^{3,\infty}} ||g^3||_{H^3} \right).
		\end{split}
		\end{equation}
		We calculate that
		\begin{equation}\label{204}
		\begin{split}
		&|| \chi^2_{M_0}||_{H^3}+|| \chi^3_{M_0}||_{H^3} \leq C  \sum_{i=0}^{3}M_0^{-i}M_0^{\frac{3}{2}} \leq C  M_0^{\frac{3}{2}},
		\\
		& ||\nabla(\chi^2_{M_0})||_{W^{3,\infty}}+||\nabla(\chi^3_{M_0})||_{W^{3,\infty}} \leq C  \sum_{i=0}^{3}M_0^{-i-1}\leq C M_0^{-1}.
		\end{split}
		\end{equation}
		For $f\times f=0, \ g \times g=0$, then we have
		\begin{align}\label{205}
		&||f \times (\nabla \times f)||_{W^{3,\infty}} + ||g \times (\nabla \times g)||_{W^{3,\infty}}
		\nonumber
		\\
		=&||f \times (\nabla \times f-f)||_{W^{3,\infty}} + ||g \times (\nabla \times g-g)||_{W^{3,\infty}}
		\nonumber
		\\
		\leq  & ||f||_{W^{3,\infty}} ||\nabla \times f-f||_{W^{3,\infty}}+||g||_{W^{3,\infty}}||\nabla \times g-g||_{W^{3,\infty}}
		\nonumber
		\\
		\leq & C \left(||(1+|\xi|)^3 \hat{f} ||_{L^1_\xi} ||(1+|\xi|)^3 (|\xi|-1)\hat{f} ||_{L^1_\xi}\right)
		\nonumber
		\\
		&\quad \quad + C||(1+|\xi|)^3 \hat{g} ||_{L^1_\xi} ||(1+|\xi|)^3 (|\xi|-1)\hat{g} ||_{L^1_\xi}
		\nonumber
		\\
		\leq &C M_1^2 \delta \left( \alpha_1^2 e^{-\frac{\nu }{2^{\alpha-1}}t}+\alpha_2^2 e^{-\frac{\mu t}{2}}\right),
		\end{align}
		and
		\begin{equation}\label{g1}
		\begin{split}
		||\left( (\nabla \times g) \times g \right) \times g ||_{W^{3,\infty}} &\leq ||\nabla \times g-g||_{W^{3,\infty}}||g||^2_{W^{3,\infty}}
		\\
		& \leq C\delta ||\hat{g}(\xi)||^3_{L^1_{\xi}}
		\\
		& \leq C\delta M_1^3 |\alpha_2|^3 e^{-\frac{3\mu t}{4}}.
		\end{split}
		\end{equation}
		Noticing that $\text{supp} \ \widehat{|f|^2}, \text{supp}\ \widehat{|g|^2} \subseteq \left\{ \xi| |\xi| \leq 2+2\delta  \right\}$, $0 < \delta \leq \frac{1}{2}$, we derive that
		\begin{align}\label{206}
		|| |f|^2 ||_{H^3} + || |g|^2 ||_{H^3}& \leq C  || |f|^2 ||_{L^2} + || |g|^2 ||_{L^2}
		\nonumber
		\\
		& \leq C  \left(||f||^2_{L^4}+||g||^2_{L^4} \right)
		\nonumber
		\\
		& \leq C \left(  \alpha_1^2 e^{-\frac{\nu t}{2^{\alpha-1}}} M_2^2 + \alpha_2^2 e^{-\frac{\mu t}{2}} M_2^2 \right),
		\end{align}
		and
		\begin{equation}\label{g0}
		\begin{split}
		|| |g|^3 ||_{H^3}& \leq C   || |g|^3 ||_{L^2}
		 \leq C  ||g||^3_{L^6}  \leq C  |\alpha_2|^3 M_2^3 e^{-\frac{3\mu t}{4}} .
		\end{split}
		\end{equation}
		Combining inequalities \eqref{204}, \eqref{205} and \eqref{206}, we get
		\begin{align*}
		& ||\tilde{f} \times \left( \nabla \times \tilde{f}\right)||_{H^3}
		+||\tilde{g} \times \left( \nabla \times \tilde{g}\right)||_{H^3}
		\leq
		C\left( \alpha_1^2 e^{-\frac{\nu t}{2^{\alpha-1}}}+\alpha_2^2 e^{-\frac{\mu t}{2}}\right)\left( \delta M_0^{\frac{3}{2}}M_1^2+M_0^{-1}M_2^2\right) .
		\end{align*}
		Combining inequalities \eqref{g2}, \eqref{g1} and \eqref{g0}, we deduce that
		\begin{equation*}
		||\left( (\nabla \times \tilde{g}) \times \tilde{g} \right) \times \tilde{g} ||_{H^3} \leq C\left( \delta M_0^{\frac{3}{2}}M_1^3+M_0^{-1}M_2^3 \right)|\alpha_2|^3 e^{-\frac{3\mu t}{4}}.
		\end{equation*}
		In what follows, we will estimate $\int^t_0 || \tilde{f} \times \tilde{g}||_{H^3}(\tau)d\tau.$
		On one hand,
		\begin{equation}\label{207}
		|| \tilde{f} \times \tilde{g}||_{H^3}= || \chi_{M_0}f \times (\chi_{M_0}g)||_{H^3} \leq   ||\chi_{M_0}^2||_{H^3} ||f \times g||_{W^{3, \infty}}.
		\end{equation}
		On the other hand, $\text{supp} \ \widehat{f \times g} \subseteq \left\{ \xi| |\xi| \leq 2+2\delta  \right\}$, $0 < \delta \leq \frac{1}{2}$. Then we have
		\begin{equation}\label{208}
		|| \tilde{f} \times \tilde{g}||_{H^3} \leq C  M_0^{\frac{3}{2}} ||\widehat{f \times g} ||_{L^1_\xi}.
		\end{equation}
		Calculate
		\begin{equation}\label{f0}
		\begin{split}
		\widehat{f \times g} & = \alpha_1 \alpha_2 \int_{\mathbb{R}^3} e^{-\nu|\xi-\eta|^\alpha t} \hat{v_0}(\xi-\eta) \times e^{-\mu|\eta|^2 t} \hat{v_0}(\eta)d\eta
		\\
		&= \frac{1}{2}\alpha_1 \alpha_2 \int_{\mathbb{R}^3} \left( e^{-(\nu|\xi-\eta|^\alpha+\mu|\eta|^2) t}-e^{-(\mu|\xi-\eta|^2+\nu|\eta|^\alpha) t} \right) \hat{v_0}(\xi-\eta) \times  \hat{v_0}(\eta)d\eta,
		\end{split}
		\end{equation}
		and
		\begin{equation}\label{209}
		\begin{split}
		& \left| e^{-(\nu|\xi-\eta|^\alpha+\mu|\eta|^2) t}-e^{-(\mu|\xi-\eta|^2+\nu|\eta|^\alpha) t} \right|\\
		=& e^{-\nu|\xi-\eta|^\alpha t} |e^{-\mu|\eta|^2t}-e^{-\mu|\xi-\eta|^2t} |+
		e^{-\mu|\xi-\eta|^2 t} |e^{-\nu|\xi-\eta|^\alpha t}-e^{-\nu|\eta|^\alpha t} |
		\\
		\leq & C t e^{-\nu |\xi-\eta|^\alpha t} \left| |\xi-\eta|^2-|\eta|^2\right|+C t e^{-\mu |\xi-\eta|^2 t} \left| |\xi-\eta|^\alpha -|\eta|^\alpha\right|\\
		\leq & C e^{-\frac{\nu}{2} |\xi-\eta|^\alpha t} \frac{\left| |\xi-\eta|^2-|\eta|^2\right|}{|\xi-\eta|^\alpha}+C e^{-\frac{\mu}{2} |\xi-\eta|^2 t} \frac{\left| |\xi-\eta|^\alpha -|\eta|^\alpha\right|}{|\xi-\eta|^2}.
		\end{split}
		\end{equation}
		\qquad In the support of  $\hat{v_0}(\xi-\eta) \times  \hat{v_0}(\eta)$, we have
		\begin{equation}\label{702}
		\frac{||\xi-\eta|^2-|\eta|^2|}{|\xi-\eta|^\alpha} \leq 3^{1-\alpha} \delta, \quad \frac{||\xi-\eta|^\alpha-|\eta|^\alpha|}{|\xi-\eta|^2} \leq 8 \alpha \delta.
		\end{equation}
		Therefore, we conclude that
		\begin{equation}\label{f1}
		 \int^\infty_0 || \tilde{f} \times \tilde{g}||_{H^3}(t)dt \leq CM_0^{\frac{3}{2}}M_1^2 (1+\alpha)\delta
		 \end{equation}
		  Then we complete the proof of Lemma \ref{tuilun2}.
	\end{proof}
	\section{The proof of Theorem \ref{thm} }
	In this section, we will prove Theorem \ref{thm} using a perturbation argument along with a standard cut-off technique.
	\begin{proof}[\\Proof of Theorem \ref{thm}]
		Let $\tilde{f}=\chi_{M_0}f, \tilde{g}=\chi_{M_0}g$, and $u=U+\tilde{f}, b=B+\tilde{g}$. Then $U,B$ satisfy
		\begin{equation}\label{U}
			\begin{split}
				&U_t+\nu \Lambda^\alpha U+\nabla\left(p+\frac{1}{2}|\tilde{f}|^2-\frac{1}{2}|\tilde{g}|^2 \right)
				\\
				=&-U \cdot \nabla U-\tilde{f} \cdot \nabla U - U \cdot \nabla \tilde{f}+B \cdot \nabla B
				\\ & \quad \quad+ \tilde{g} \cdot \nabla B+B \cdot \nabla \tilde{g}+F,
			\end{split}
		\end{equation}
		
		\begin{equation}\label{B}
			\begin{split}
				&B_t-\mu \Delta B- \nabla \times (((\nabla \times B) \times B) \times B)
				\\
				=&-U \cdot \nabla B-\tilde{f} \cdot \nabla B - U \cdot \nabla \tilde{g}+B \cdot \nabla U+ \tilde{g} \cdot \nabla U+ B \cdot \nabla \tilde{f}
				\\ & \quad -\eta \nabla \times \left( \left(\nabla \times B \right) \times B\right)-\eta \nabla \times \left( \left(\nabla \times B \right) \times \tilde{g}\right)
				\\
				& \quad - \eta \nabla \times \left( \left(\nabla \times \tilde{g} \right) \times B\right)-\eta \nabla \times \left( (\nabla \times \tilde{g})\times \tilde{g}\right)+G
				\\
				&\quad +\nabla \times (((\nabla \times B) \times B) \times \tilde{g})+\nabla \times (((\nabla \times B) \times \tilde{g}) \times B)
				\\
				& \quad +\nabla \times (((\nabla \times B) \times \tilde{g}) \times \tilde{g}) +\nabla \times (((\nabla \times \tilde{g}) \times B) \times B)
				\\
				& \quad +\nabla \times (((\nabla \times \tilde{g}) \times B) \times \tilde{g}) +\nabla \times (((\nabla \times \tilde{g}) \times \tilde{g}) \times B),
			\end{split}
		\end{equation}
		where
		\begin{equation}\label{fg}
		\begin{split}
		&F:= \tilde{f}\times (\nabla \times \tilde{f})-\tilde{g}\times \left(\nabla \times \tilde{g}\right)-\nu \Delta\chi_{M_0}f+2\nu \nabla \cdot
		\left( \nabla\chi_{M_0} f\right),
		\\
		&G:= \nabla \times ( \tilde{f} \times \tilde{g})-\mu \Delta\chi_{M_0}g+2\mu \nabla \cdot
		\left( \nabla\chi_{M_0} g\right)+\frac{1}{2}f \cdot \nabla {\chi_{M_0}}^2g
		\\
		& \quad \quad \quad-\frac{1}{2}g \cdot \nabla {\chi_{M_0}}^2f-\nabla \times (((\nabla \times \tilde{g}) \times \tilde{g} ) - \nabla \times (((\nabla \times \tilde{g}) \times \tilde{g} ) \times \tilde{g}).
		\end{split}
		\end{equation}
		In what follows, we will derive some energy estimates of $U$ and $B$.

		\textbf{\\Step 1: Energy inequalities of $B$.\\}
		Taking the derivatives $\Lambda^{k}, 0 \leq k \leq 3$ on Equation \ref{B} and $L^2$ inner product with $\Lambda^{k} B$, we get
		\begin{equation}\label{307}
		\begin{split}
			&\frac{1}{2}\frac{d}{dt}|| B||^2_{H^3}+ || \nabla   B||^2_{H^3}
			\nonumber
			\\
			= & I_1+I_2+I_3+I_4+I_5+I_6
			\\
			&\quad+ J_1+J_2+J_3 + \int_{\mathbb{R}^3} G \cdot B dx
			\\
			& \quad+ K_0+K_1+K_2+K_3+K_4+K_5+K_6,
			\end{split}
		\end{equation}
		where
		\begin{equation*}
			\begin{split}
				&I_1=-\sum_{0 \leq k\leq 3}\int_{\mathbb{R}^3} \Lambda^k (U \cdot \nabla B) \cdot \Lambda^k B dx, \quad I_2=-\sum_{0 \leq k\leq 3}\int_{\mathbb{R}^3} \Lambda^k(\tilde{f} \cdot \nabla B) \cdot \Lambda^k B dx,
				\\
				&I_3=-\sum_{0 \leq k\leq 3} \int_{\mathbb{R}^3} \Lambda^k (U \cdot \nabla \tilde{g}) \cdot \Lambda^k B dx, \quad
				I_4=\sum_{0 \leq k\leq 3} \int_{\mathbb{R}^3} \Lambda^k (B \cdot \nabla U) \cdot \Lambda^k B dx, \\
				& I_5=\sum_{0 \leq k\leq 3}  \int_{\mathbb{R}^3} \Lambda^k (\tilde{g} \cdot \nabla U) \cdot \Lambda^k B dx, \quad I_6=\sum_{0 \leq k\leq 3} \int_{\mathbb{R}^3} \Lambda^k (B \cdot \nabla \tilde{f}) \cdot \Lambda^k B dx,
				\end{split}
				\end{equation*}
			\begin{equation*}
			\begin{split}
				&J_1= -\sigma \sum_{0 \leq k\leq 3} \int_{\mathbb{R}^3}  \Lambda^k  (( \nabla \times B) \times B) \cdot (\nabla \times \Lambda^k B) dx, \\
				& J_2=-\sigma \sum_{0 \leq k\leq 3} \int_{\mathbb{R}^3}  \Lambda^k (( \nabla \times B) \times \tilde{g}) \cdot (\nabla \times \Lambda^k B) dx,
				\\
				& J_3=-\sigma \sum_{0 \leq k\leq 3} \int_{\mathbb{R}^3}  \Lambda^k(( \nabla \times \tilde{g}) \times B) \cdot (\nabla \times \Lambda^k B) dx,
			\end{split}
		\end{equation*}
\begin{equation}\label{22}
\begin{split}
&K_0=\kappa \sum_{0 \leq k\leq 3}\int_{\mathbb{R}^3}  \Lambda^k(((\nabla \times B) \times B) \times B)\cdot (\nabla \times \Lambda^k B) dx,
\\
& K_1=\kappa \sum_{0 \leq k\leq 3}\int_{\mathbb{R}^3} \Lambda^k(((\nabla \times B) \times B) \times \tilde{g})\cdot (\nabla \times \Lambda^k B) dx,
\\
& K_2=\kappa \sum_{0 \leq k\leq 3} \int_{\mathbb{R}^3} \Lambda^k(((\nabla \times B) \times \tilde{g}) \times B) \cdot (\nabla \times \Lambda^k B) dx,
\\
& K_3=\kappa \sum_{0 \leq k\leq 3} \int_{\mathbb{R}^3} \Lambda^k(((\nabla \times B) \times \tilde{g}) \times \tilde{g})\cdot (\nabla \times \Lambda^k B) dx,
\end{split}
\end{equation}
\begin{equation*}
\begin{split}
& K_4=\kappa \sum_{0 \leq k\leq 3} \int_{\mathbb{R}^3} \Lambda^k(((\nabla \times \tilde{g}) \times B) \times \tilde{g}) \cdot (\nabla \times \Lambda^k B) dx,
\\
& K_5=\kappa \sum_{0 \leq k\leq 3} \int_{\mathbb{R}^3} \Lambda^k(((\nabla \times \tilde{g}) \times B) \times B)\cdot (\nabla \times \Lambda^k B) dx,
\\
&K_6=\kappa \sum_{0 \leq k\leq 3} \int_{\mathbb{R}^3} \Lambda^k(((\nabla \times \tilde{g}) \times \tilde{g}) \times B) \cdot (\nabla \times \Lambda^k B) dx.
\end{split}
\end{equation*}
Firstly, we estimate $I_1,I_2$ in the following:
\begin{equation*}
\begin{split}
&|I_1+I_2|
\\
\leq &\sum_{0 \leq k\leq 3} \left|\left(\int_{\mathbb{R}^3} (\Lambda^k (U \cdot \nabla B)- (U \cdot \nabla \Lambda^k B) \cdot \Lambda^k B dx +\int_{\mathbb{R}^3} ( \Lambda^k(\tilde{f} \cdot \nabla B- \tilde{f} \cdot \nabla \Lambda^k B)) \cdot \Lambda^k B dx \right)  \right|
\\
& \quad + \big|\sum_{0 \leq k\leq 3} \int_{\mathbb{R}^3}  (u \cdot \nabla \Lambda^k B) \cdot \Lambda^k B dx \big|
\\
\leq &\sum_{0 \leq k\leq 3} \left(\big|\int_{\mathbb{R}^3} (\Lambda^k (U \cdot \nabla B)- (U \cdot \nabla \Lambda^k B) \cdot \Lambda^k B dx \big| + \big|\int_{\mathbb{R}^3} ( \Lambda^k(\tilde{f} \cdot \nabla B- \tilde{f} \cdot \nabla \Lambda^k B)) \cdot \Lambda^k B dx \big| \right)
\\
 \leq &C \left(||\nabla U||_{L^\infty}||\nabla B||_{H^2}+||\nabla B||_{L^\frac{6}{\alpha}} ||U||_{W^{3,\frac{6}{2-\alpha}}}\right) || B||_{H^3}
 \\
 & \quad+ C \left( ||\nabla B||_{H^2} ||\nabla \tilde{f} ||_{L^\infty}+ ||\nabla B||_{L^6} ||\tilde{f}||_{W^{3,3}}\right)|| B||_{H^3}.
\end{split}
\end{equation*}
By Sobolev's inequality, we deduce that
\begin{equation}\label{I12}
|I_1+I_2| \leq C \left(||\Lambda^\frac{\alpha}{2}U ||_{H^3}||\nabla B ||_{H^3}|| B ||_{H^3}+||B||^2_{H^3}( ||\tilde{f}||_{W^{1,\infty}}+ ||\tilde{f}||_{W^{3,3}}) \right).
\end{equation}
For $I_3$, it's easy for us to get
\begin{equation}\label{I3}
|I_3| \leq C ||U||_{H^3}||B||_{H^3}|| \tilde{g} ||_{W^{4,\infty}}.
\end{equation}
For $I_3, I_4$, we derive that
\begin{equation}\label{I45}
\begin{split}
 |I_4+I_5| \leq &C\left( ||\nabla U||_{H^2}||\nabla B||_{L^\infty}+ ||\nabla U||_{L^3} ||B||_{W^{3,6}}\right)||B||_{H^3}
\\
& \quad + C(||\nabla U||_{H^2}||\nabla \tilde{f} ||_{L^\infty}+ ||\nabla U||_{L^6} ||\tilde{f}||_{W^{3,3}})||B||_{H^3}
\\
& \leq C \left(||\Lambda^\frac{\alpha}{2}U||_{H^3}||\nabla B||_{H^3}||B||_{H^3} + ||U||_{H^3}||B||_{H^3}( ||\tilde{f}||_{W^{1,\infty}}+ ||\tilde{f}||_{W^{3,3}}) \right).
\end{split}
\end{equation}
Considering $I_6$, we have
\begin{equation}\label{I6}
|I_6| \leq C ||B||^2_{H^3} ||\tilde{f}||_{W^{4,\infty}}.
\end{equation}
Next step, we will estimate $J_1, J_2, J_3$. For
\begin{align*}
J_1&=
 \sum_{0 \leq k\leq 3} \sigma \int_{\mathbb{R}^3} \Lambda^k (\nabla \times B) \cdot \Lambda^k  \left( \left(\nabla \times B \right) \times B\right) dx
\\
&= \sigma \sum_{0 \leq k\leq 3} \int_{\mathbb{R}^3} \left\{   \Lambda^k  \left( \left(\nabla \times B \right) \times B\right)  -   \Lambda^k (\nabla \times B ) \times B  \right\}\cdot  \Lambda^k (\nabla \times B) dx.
\end{align*}
Using Lemma \ref{lemc}, we deduce that
\begin{align}\label{J1}
|J_1| &\leq C\eta \sum_{0 \leq k\leq 3} ||\Lambda^k (\nabla \times B)||_{L^2} ||\Lambda^k  \left( \left(\nabla \times B \right) \times B\right)  -  \Lambda^k (\nabla \times B ) \times B   ||_{L^2}
\nonumber
\\
&\leq C \eta ||\nabla B||_{H^3}\left( ||\nabla \times B||_{H^2}||\nabla B||_{L^\infty}+||\nabla \times B||_{L^\infty}||B||_{H^3}\right)
\nonumber
\\
&\leq C \eta||\nabla B||^2_{H^3}|| B||_{H^3}.
\end{align}
For $J_2$, we calculate
\begin{align*}
J_2
&=\sum_{0 \leq k\leq 3} \sigma \int_{\mathbb{R}^3} \Lambda^k (\nabla \times B) \cdot \Lambda^k  \left( \left(\nabla \times B \right) \times \tilde{g}\right)dx
\\
&=\sigma \sum_{0 \leq k\leq 3}  \int_{\mathbb{R}^3} \Lambda^k (\nabla \times B) \cdot \left( \Lambda^k  \left( \left(\nabla \times B \right) \times \tilde{g}\right) - \Lambda^k(\nabla \times B ) \times \tilde{g} \right) dx
\\
& \quad \quad+\sigma \sum_{0 \leq k\leq 3} \int_{\mathbb{R}^3}  (\nabla \times \Lambda^k B) \cdot \left( (\nabla \times \Lambda^k B) \times  \tilde{g}\right)dx
\\
&=\sigma \sum_{0 \leq k\leq 3}  \int_{\mathbb{R}^3} \Lambda^k (\nabla \times B) \cdot \left( \Lambda^k  \left( \left(\nabla \times B \right) \times \tilde{g}\right) - \Lambda^k(\nabla \times B ) \times \tilde{g} \right) dx,
\end{align*}
for we use the fact that$ \int_{\mathbb{R}^3}  (\nabla \times \Lambda^k B) \cdot \left( (\nabla \times \Lambda^k B) \times  \tilde{g}\right)dx=0.$ Using Lemma \ref{lemc}, we then derive that
\begin{equation}\label{J2}
\begin{split}
|J_2| &\leq C \eta||\nabla B||_{H^3}||\nabla B||_{H^2}||\tilde{g}||_{W^{3,\infty}} \\
&  \leq C \eta||\nabla B||_{H^3} || B||_{H^3} ||\tilde{g}||_{W^{3,\infty}}
\\
&  \leq  \frac{\mu}{16}||\nabla B||^2_{H^3}+C|| B||^2_{H^3} ||\tilde{g}||^2_{W^{3,\infty}}.
\end{split}
\end{equation}
For $J_3$, we could estimate it directly that
\begin{equation}\label{J3}
\begin{split}
|J_3| & \leq C ||\tilde{g} ||_{W^{4,\infty}} || B||_{H^3} || \nabla B||_{H^3}
\\
&  \leq  \frac{\mu}{16}||\nabla B||^2_{H^3}+C|| B||^2_{H^3} ||\tilde{g}||^2_{W^{4,\infty}}.
\end{split}
\end{equation}
	For some quadratic $\nabla \Lambda^3B$(highest derivatives) in $K_1, K_2$, we could not get good estimate of $K_1, K_2$ when $\tilde{g}$ is large.  We also should not neglect $K_0, K_3$ containing some positive items. Therefore, we then find that it's a effective way to estimate $K_1, K_2, K_0$ and $K_3$ together.
		\begin{equation}\label{21}
		\begin{split}
		&K_0+K_1+K_2+K_3
		\\
		=&-\kappa\sum_{0 \leq k\leq 3}\int_{\mathbb{R}^3}  |(\nabla \times \Lambda^k B) \times B|^2+2 ((\nabla \times \Lambda^k B) \times B) \cdot (\tilde{g} \times (\nabla \times \Lambda^k B) )
		\\
		& \quad  +F_0+F_1+F_2+F_3
		-(\nabla \times \Lambda^k B) \times \tilde{g}|^2 dx
		\\
		\leq & -\kappa \sum_{0 \leq k\leq 3} \int_{\mathbb{R}^3} \left( |(\nabla \times \Lambda^\alpha B) \times B|-|(\nabla \times \Lambda^\alpha B) \times \tilde{g}| \right)^2+F_0+F_1+F_2+F_3,
		\end{split}
		\end{equation}
	where \begin{equation}\label{F0}
	F_0=\kappa\sum_{0 \leq k\leq 3} \sum_{|\beta| \leq 2, \beta + \gamma + \theta =k}
	\int_{\mathbb{R}^3} (((\nabla \times \Lambda^\beta B) \times \Lambda^ \gamma B) \times \Lambda^\theta B) \cdot (\nabla \times \Lambda^k B) dx,
	\end{equation}
	\begin{equation}\label{F1}
	F_1=\kappa\sum_{0 \leq k\leq 3} \sum_{|\beta| \leq 2, \beta + \gamma + \theta =k}
	\int_{\mathbb{R}^3} (((\nabla \times \Lambda^\beta B) \times \Lambda^ \gamma B) \times \Lambda^\theta \tilde{g}) \cdot (\nabla \times \Lambda^k B) dx,
	\end{equation}
	\begin{equation}\label{F2}
	F_2=\kappa\sum_{0 \leq k\leq 3} \sum_{|\beta| \leq 2, \beta + \gamma + \theta =k}
	\int_{\mathbb{R}^3} (((\nabla \times \Lambda^\beta B) \times \Lambda^ \gamma \tilde{g}) \times \Lambda^\theta B) \cdot (\nabla \times \Lambda^k B) dx,
	\end{equation}
	\begin{equation}\label{F3}
	F_3=\kappa\sum_{0 \leq k\leq 3} \sum_{|\beta| \leq 2, \beta + \gamma + \theta =k}
	\int_{\mathbb{R}^3} (((\nabla \times \Lambda^\beta B) \times \Lambda^ \gamma \tilde{g}) \times \Lambda^\theta \tilde{g}) \cdot (\nabla \times \Lambda^k B) dx.
	\end{equation}
	For $F_0$, we have
		\begin{equation}\label{23}
		\begin{split}
		F_0 &\leq C \left( ||\nabla B||_{H^2} || \nabla B||_{L^\infty} ||B||_{L^\infty} + ||\nabla B||_{L^\infty} || \nabla B||_{H^2} || B||_{L^\infty} + ||\nabla B||^3_{L^6} \right)||\nabla B||_{H^3}
		\\
		& \leq C||\nabla B||^2_{H^3}||B||^2_{H^3}.
		\end{split}
		\end{equation}
		For $F_1$, we derive that
		\begin{equation}\label{24}
		\begin{split}
		F_1 &\leq C ||\tilde{g}||_{W^{4,\infty}} ||B||^2_{H^3} ||\nabla B||_{H^3}
		\\
		& \leq \frac{\mu}{16}||\nabla B||^2_{H^3}+C||\tilde{g}||^2_{W^{4,\infty}} ||B||^4_{H^3}.
		\end{split}
		\end{equation}
		For $F_2$, we could get
		\begin{equation}\label{25}
		\begin{split}
		F_2 &\leq C ||\nabla B||_{H^2}\left( ||\nabla B||_{L^\infty} || \tilde{g} ||_{L^\infty} + ||B||_{L^\infty} || \tilde{g} ||_{W^{1,\infty}}\right) ||\nabla B||_{H^3}
		\\
		& \qquad + C ||\nabla B||_{W^{1,6}}\left( ||\nabla B||_{W^{1,6}} ||\tilde{g}||_{W^{1,6}}+ || B||_{W^{1,6}} ||\tilde{g}||_{W^{2,6}}\right)||\nabla B||_{H^3}
		\\
		& \leq \frac{\mu}{16}||\nabla B||^2_{H^3}+C\left( ||\tilde{g}||^2_{W^{2,\infty}}+||\tilde{g}||^2_{W^{2,6}} \right) ||B||^4_{H^3}.
		\end{split}
		\end{equation}
		For $F_3$, it's easy for us th get
		\begin{equation}\label{26}
		\begin{split}
		F_3 &\leq C  ||\nabla B||_{H^2} ||\tilde{g}||^2_{L^\infty}||\nabla B||_{H^3} +||\nabla B||_{W^{1,6}} ||\tilde{g}||^2_{W^{2,6}} ||\nabla B||_{H^3}
		\\
		& \qquad +C || \nabla B||_{L^2} ||\tilde{g}||^2_{W^{3,\infty}} ||\nabla B||_{H^3}
		\\
		& \leq \frac{\mu}{16}||\nabla B||^2_{H^3}+C \left(||\tilde{g}||^4_{W^{3,\infty}}+||\tilde{g}||^4_{W^{2,6}}  \right) ||B||^2_{H^3}.
		\end{split}
		\end{equation}
		To estimate $K_4,K_5,K_6$, we have
		\begin{equation}\label{27}
		\begin{split}
		K_4 &\leq C ||\tilde{g}||_{W^{4,\infty}} ||B||_{H^3} ||\tilde{g}||_{W^{3,\infty}} ||\nabla B||_{H^3}
		\\
		& \leq \frac{\mu}{16}||\nabla B||^2_{H^3}+C||\tilde{g}||^4_{W^{4,\infty}} ||B||^2_{H^3},
		\end{split}
		\end{equation}
		\begin{equation}\label{28}
		\begin{split}
		K_5 &\leq C ||\tilde{g}||_{W^{4,\infty}} ||B||^2_{H^3} ||\nabla B||_{H^3}
		 \\
		 & \leq \frac{\mu}{16}||\nabla B||^2_{H^3}+C||\tilde{g}||^2_{W^{4,\infty}} ||B||^4_{H^3},
		\end{split}
		\end{equation}
		\begin{equation}\label{29}
		\begin{split}
		K_6 &\leq C ||(\nabla \times \tilde{g}) \times \tilde{g}||_{W^{3,\infty}} ||B||_{H^3} ||\nabla B||_{H^3}
		\\
		& \leq \frac{\mu}{16}||\nabla B||^2_{H^3}+C||\tilde{g}||^4_{W^{4,\infty}} ||B||^2_{H^3}.
		\end{split}
		\end{equation}
		At last, we consider the term $\int_{\mathbb{R}^3}\Lambda^k G \Lambda^k B dx$. Recalling the expression  of $G$, we have
		\begin{align}\label{3007}
			\sum_{0 \leq k\leq 3} \left| \int_{\mathbb{R}^3}\Lambda^k G \Lambda^k B dx \right| &=\sum_{0 \leq k\leq 3}\big| \int_{\mathbb{R}^3} \Lambda^k (\nabla \times ( \tilde{f} \times \tilde{g})+2\nu\nabla \cdot \left( \nabla\chi_{M_0}g\right)-\nu\Delta \chi_{M_0}g
			\nonumber
			\\
			& \quad \quad+\frac{1}{2}f \cdot \nabla \chi^2_{M_0}g-\frac{1}{2}g \cdot \nabla \chi^2_{M_0}f ) \cdot \Lambda^k B dx\big|
			\nonumber
			\\
			& \quad \quad+ \sum_{0 \leq k\leq 3}|\int_{\mathbb{R}^3} \Lambda^k \left(  \sigma(\nabla \times \tilde{g}) \times \tilde{g})-\kappa((\nabla \times \tilde{g}) \times \tilde{g}) \times \tilde{g} ) \right) \cdot \Lambda^k (\nabla \times B) dx|
			\nonumber
			\\
			&  \leq C ||\tilde{f}\times \tilde{g}||_{H^3}||\nabla B||_{H^3}+||\nabla\chi_{M_0}g||_{H^3}||\nabla B||_{H^3}
			+C||\Delta \chi_{M_0}g||_{W^{3,\frac{6}{5}}}||B||_{W^{3,6}}
			\nonumber
			\\
			& \quad \quad+C\left(||f \cdot \nabla \chi^2_{M_0}g||_{W^{3,\frac{6}{5}}}+||g \cdot \nabla \chi^2_{M_0}f||_{W^{3,\frac{6}{5}}}\right)||B||_{W^{3,6}}
			\nonumber
			\\
			& \quad \quad +C\left(||(\nabla \times \tilde{g}) \times \tilde{g}||_{H^3}+C||((\nabla \times \tilde{g}) \times \tilde{g})\times \tilde{g}||_{H^3} \right) ||\nabla B||_{H^3}
			\nonumber
			\\
			& \leq C \left(||\tilde{f}\times \tilde{g}||_{H^3}+M_0^{-1}||g||_{H^3(|x|\leq  2M_0)}  \right)||\nabla B||_{H^{3}}
			\nonumber
			\\
			& \quad \quad + C M_0^{-2}||g||_{W^{3,\frac{6}{5}}(|x|\leq 2M_0)}  ||\nabla B||_{H^{3}}
			\nonumber
			\\
			& \quad \quad + CM^{-1}_0 ||f\otimes g||_{W^{3,\frac{6}{5}}(|x|\leq 2M_0)}  ||\nabla B||_{H^{3}}
			\nonumber
			\\
			& \leq C \left( ||\tilde{f}\times \tilde{g}||_{H^3}+ M_0^{-\frac{1}{2}}M_2 |\alpha_2| e^{-\frac{\mu t}{4}}\right) ||\nabla B||_{H^{3}}
			\nonumber
			\\
			& \quad \quad+ C|\alpha_1\alpha_2| M_0^{-\frac{1}{2}} M_2^2 e^{-\frac{(\mu
					+\nu)t}{4}}  ||\nabla B||_{H^{3}}
			\\
			& \quad \quad+C\left(||(\nabla \times \tilde{g}) \times \tilde{g}||_{H^3}+C||((\nabla \times \tilde{g}) \times \tilde{g})\times \tilde{g}||_{H^3} \right) ||\nabla B||_{H^3}.
			\nonumber
		\end{align}

		\textbf{\\Step 2: Energy inequalities of $U$.\\}
		Operating Equation \ref{U} with $\Lambda^k, 0 \leq k\leq 3$, and taking $L^2$ on Equation \ref{U} yields
		\begin{equation}\label{314}
			\begin{split}
				&\frac{1}{2}\frac{d}{dt}|| U||^2_{H^3}+\mu ||\nabla U||^2_{H^3}
				\nonumber
				\\
				= & H_1+H_2+H_3+H_4+H_5+H_6
				\\
				& \quad+ \sum_{0 \leq k\leq 3}\int_{\mathbb{R}^3} \Lambda^k F \Lambda^k Udx-\sum_{0 \leq k\leq 3} \int_{\mathbb{R}^3}\Lambda^k \nabla \left(p+\frac{1}{2}|\tilde{f}|^2-\frac{1}{2}|\tilde{g}|^2 \right)\Lambda^k Udx,
			\end{split}
		\end{equation}
where
\begin{equation*}
\begin{split}
&H_1=-\sum_{0 \leq k\leq 3}\int_{\mathbb{R}^3} \Lambda^k (U \cdot \nabla U) \cdot \Lambda^k U dx, \quad H_2=-\sum_{0 \leq k\leq 3}\int_{\mathbb{R}^3} \Lambda^k(\tilde{f} \cdot \nabla U) \cdot \Lambda^k U dx,
\\
&H_3=-\sum_{0 \leq k\leq 3} \int_{\mathbb{R}^3} \Lambda^k (U \cdot \nabla \tilde{f}) \cdot \Lambda^k U dx, \quad
H_4=\sum_{0 \leq k\leq 3} \int_{\mathbb{R}^3} \Lambda^k (B \cdot \nabla B) \cdot \Lambda^k U dx, \\
& H_5=\sum_{0 \leq k\leq 3}  \int_{\mathbb{R}^3} \Lambda^k (\tilde{g} \cdot \nabla B) \cdot \Lambda^k U dx, \quad H_6=\sum_{0 \leq k\leq 3} \int_{\mathbb{R}^3} \Lambda^k (B \cdot \nabla \tilde{g}) \cdot \Lambda^k U dx,
\end{split}
\end{equation*}
		Firstly, we have
\begin{equation*}
\begin{split}
&|H_1+H_2|
\\
\leq &\left|\sum_{0 \leq k\leq 3} \left(\int_{\mathbb{R}^3} (\Lambda^k (U \cdot \nabla U)- (U \cdot \nabla \Lambda^k U) \cdot \Lambda^k U dx +\int_{\mathbb{R}^3} ( \Lambda^k(\tilde{f} \cdot \nabla U- \tilde{f} \cdot \nabla \Lambda^k U)) \cdot \Lambda^k U dx \right)  \right|
\\
& \quad + \left|\sum_{0 \leq k\leq 3} \int_{\mathbb{R}^3}  \left( (U \cdot \nabla \Lambda^k U) \cdot \Lambda^k U + (\tilde{f} \cdot \nabla \Lambda^k U)\cdot \Lambda^k U \right)dx \right|
\\
\leq &\sum_{0 \leq k\leq 3} \left|\left(\int_{\mathbb{R}^3} (\Lambda^k (U \cdot \nabla U)- (U \cdot \nabla \Lambda^k U) \cdot \Lambda^k U dx +\int_{\mathbb{R}^3} ( \Lambda^k(\tilde{f} \cdot \nabla U- \tilde{f} \cdot \nabla \Lambda^k U)) \cdot \Lambda^k U dx \right)  \right|
\\
& \quad + \big|\sum_{0 \leq k\leq 3} \int_{\mathbb{R}^3}  (u \cdot \nabla \Lambda^k U) \cdot \Lambda^k U dx \big|
\\
\leq &\big|\sum_{0 \leq k\leq 3} \left(\int_{\mathbb{R}^3} (\Lambda^k (U \cdot \nabla U)- (U \cdot \nabla \Lambda^k U) \cdot \Lambda^k U dx \big| + \big|\int_{\mathbb{R}^3} ( \Lambda^k(\tilde{f} \cdot \nabla U- \tilde{f} \cdot \nabla \Lambda^k U)) \cdot \Lambda^k U dx \big| \right)
\\
\leq &C \left(||\nabla U||_{L^\infty}||\nabla U||_{H^2}+||\nabla U||_{L^\frac{6}{\alpha}} ||U||_{W^{3,\frac{6}{2-\alpha}}}\right) || U||_{H^3}
\\
& \quad+ C \left( ||\nabla U||_{H^2} ||\nabla \tilde{f} ||_{L^\infty}+ ||\nabla U||_{L^6} ||\tilde{f}||_{W^{3,3}}\right)|| U||_{H^3},
\end{split}
\end{equation*}
By Sobolev inequality, we deduce that
\begin{equation}\label{H12}
|H_1+H_2| \leq C \left(||\Lambda^\frac{\alpha}{2}U ||^2_{H^3}||U ||_{H^3}+||U||^2_{H^3}( ||\tilde{f}||_{W^{1,\infty}}+ ||\tilde{f}||_{W^{3,3}}) \right).
\end{equation}
For $H_3$, it's easy for us to get
\begin{equation}\label{H3}
|H_3| \leq C ||U||_{H^3}||B||_{H^3}|| \tilde{f} ||_{W^{4,\infty}}.
\end{equation}
Using the similar way to estimate $H_1+H_2$, we have
\begin{equation}\label{H45}
\begin{split}
|H_4+H_5| \leq &C\left( ||\nabla B||_{H^2}||\nabla B||_{L^\infty}+ ||\nabla B||_{L^3} ||B||_{W^{3,6}}\right)||U||_{H^3}
\\
& \quad + C(||\nabla B||_{H^2}||\nabla \tilde{g} ||_{L^\infty}+ ||\nabla B||_{L^6} ||\tilde{g}||_{W^{3,3}})||U||_{H^3}
\\
& \leq C \left(||\nabla B||^2_{H^3}||B||_{H^3} + ||U||_{H^3}||B||_{H^3}( ||\tilde{g}||_{W^{1,\infty}}+ ||\tilde{g}||_{W^{3,3}}) \right)
\end{split}
\end{equation}
For $H_6$, we have
\begin{equation}\label{H6}
|H_6| \leq C ||U||_{H^3}||B||_{H^3} ||\tilde{g}||_{W^{4,\infty}}.
\end{equation}
As for $\int_{\mathbb{R}^3}\Lambda^k F \Lambda^k U dx$, it suffices for us to have
		\begin{equation}\label{30000}
			\begin{split}
				\sum_{0 \leq k\leq 3}\left| \int_{\mathbb{R}^3}\Lambda^k F \Lambda^k U dx \right| &= \sum_{0 \leq k\leq 3} \big| \int_{\mathbb{R}^3} \Lambda^k (\tilde{f}\times (\nabla \times \tilde{f})-\tilde{g}\times \left(\nabla \times \tilde{g}\right)-\nu \Delta\chi_{M_0}f
				\\
				& \quad \quad+2\nu \nabla \cdot
				\left( \nabla\chi_{M_0} f\right) ) \cdot \Lambda^k U dx\big|
				\\
				&  \leq C \left( ||\tilde{f}\times (\nabla \times \tilde{f})||_{H^3}+||\tilde{g}\times (\nabla \times \tilde{g})||_{H^3} \right)|| U||_{H^3}
				\\
				& \quad \quad+C\left(|| \nabla\chi_{M_0}f ||_{W^{4,\frac{6}{4+\alpha}}}
				+||\Delta \chi_{M_0}f||_{W^{3,\frac{6}{4+\alpha}}}\right)||U||_{W^{3,\frac{6}{2-\alpha}}}
				\\
				& \leq C \left( ||\tilde{f}\times (\nabla \times \tilde{f})||_{H^3}+||\tilde{g}\times (\nabla \times \tilde{g})||_{H^3} \right)|| U||_{H^3}
				\\
				&\quad \quad +CM_0^{-\frac{1}{2}}M_2 |\alpha_1| e^{-\frac{\nu t}{2^{\alpha}}}  ||\Lambda^{\frac{\alpha}{2}} U||_{H^{3}}.
			\end{split}
		\end{equation}
		We could estimate the pressure term in the same way with in \cite{Z}. Since
		\begin{align*}
			p&=(-\Delta)^{-1} \text{div} \left( u \cdot \nabla u-b \cdot \nabla b \right)
			\\
			&=\sum_{i,j}(-\Delta)^{-1} \partial_i \partial_j \left(u_iU_j-b_iB_j \right)+(-\Delta)^{-1} \nabla \cdot \left(U \cdot \nabla \tilde{f}-B\cdot \nabla \tilde{g} \right)
			\\
			& \quad \quad + (-\Delta)^{-1} \nabla \cdot \left( \tilde{f}\times (\nabla \times \tilde{f})-\tilde{g}\times \left(\nabla \times \tilde{g}\right) \right)-\frac{1}{2}|\tilde{f}|^2+\frac{1}{2}|\tilde{g}|^2,
		\end{align*}
		then we have
		\begin{equation*}
			\begin{split}
				\Pi&:=\left|-\sum_{|k|\leq 3}\int_{\mathbb{R}^3}\Lambda^k \nabla \left(p+\frac{1}{2}|\tilde{f}|^2-\frac{1}{2}|\tilde{g}|^2 \right)\Lambda^k Udx \right|
				\\
				&  \leq \left( || u \otimes U||_{W^{3,\frac{3}{2}}} +|| h \otimes B||_{W^{3,\frac{3}{2}}}\right) ||f \cdot \nabla\chi_{M_0}||_{W^{3,3}}
				\\
				&  \quad \quad
				+ ||U\cdot \nabla \tilde{f}-B\cdot \nabla \tilde{g}||_{H^3}||U||_{H^3}
				\\
				&  \quad \quad
				+  ||\tilde{f}\times (\nabla \times \tilde{f})-\tilde{g}\times \left(\nabla \times \tilde{g}\right)||_{H^3}||U||_{H^3}.
			\end{split}
		\end{equation*}	
		By H\"older's inequality, we furthermore derive that
		\begin{equation*}
			\begin{split}
				\Pi & \leq C \big[\left( ||U||^2_{H^3}+||\tilde{f}||_{W^{3,6}}||U||_{H^3} \right)||\nabla \chi_{M_0}||_{W^{3,\infty}} ||f||_{W^{3,3}(M_0 \leq |x| \leq 2M_0)}
				\\
				&
				\quad \quad+||H||^2_{H^3}\cdot||\nabla \chi_{M_0}||_{W^{3,\infty}} ||f||_{W^{3,3}(M_0 \leq |x| \leq 2M_0)}
				\\
				&   \quad \quad+||\tilde{g}||_{W^{3,2}} ||H||_{W^{3,6}}||\nabla \chi_{M_0}||_{W^{3,\infty}} ||f||_{W^{3,3}(M_0 \leq |x|\leq 2M_0)}
				\\
				&   \quad \quad+ \left( ||\nabla \tilde{f}||_{W^{3,\infty}}||U||_{H^3}+||\nabla \tilde{g}||_{W^{3,\infty}}||B||_{H^3}\right)
				\\
				&  \quad \quad+\left( ||\tilde{f}\times (\nabla \times \tilde{f})||_{H^3}+||\tilde{g}\times \left(\nabla \times \tilde{g}\right)||_{H^3} \right)||U||_{H^3}
				\\
				& \quad \quad +||\nabla B||_{H^3}  ||\left(\nabla \times \tilde{g} \right) \times \tilde{g}||_{H^3}+||B||^2_{H^3}||\tilde{g}||_{W^{5,\infty}} \big],
				\\
				&  \leq C\big[ ( ||U||^2_{H^3} +|\alpha_1| M_2 e^{-\frac{\mu t}{4}}||U||_{H^3} + ||B||^2_{H^3}) \alpha_1 M_2 M_0^{-1}e^{-\frac{\mu t}{4}}
				\\
				& \quad \quad
				+ |\alpha_2| M_2 M_0^{\frac{1}{2}} e^{ -\frac{\nu t}{2^{\alpha}}} ||\nabla B||_{H^3}  \cdot \alpha_1 M_2 M_0^{-1}e^{-\frac{\mu t}{4}}
				\nonumber
				\\
				&   \quad \quad+ \left( |\alpha_1|M_1 e^{-\frac{\mu t}{4}}||U||_{H^3}+|\alpha_2| M_1 e^{-\frac{\nu t}{2^{\alpha}}}||B||_{H^3}\right)||U||_{H^3}
				\nonumber
				\\
				&   \quad \quad +\left( \alpha_1^2 e^{-\frac{\nu t}{2^{\alpha-1}}}+\alpha_2^2 e^{-\frac{\mu t}{2}}\right) \cdot \left( \delta M_0^{\frac{3}{2}}M_1^2+M_0^{-1}M_2^2 \right) ||U||_{H^3}\big] .
			\end{split}
		\end{equation*}
		In a result, we get
		\begin{equation}\label{315}
			\begin{split}
				\Pi &\leq C \left(|\alpha_1| (M_1+M_2)e^{-\frac{\mu t}{4}} ||U||^2_{H^3}+|\alpha_1| M_2 e^{-\frac{\mu t}{4}} ||B||^2_{H^3}\right)
				\\
				& \quad \quad+
				C|\alpha_2| M_1 e^{-\frac{\nu t}{2^{\alpha}}}||U||_{H^3}||B||_{H^3}
				\\
				& \quad \quad+C\left( \alpha_1^2 e^{-\frac{\nu t}{2^{\alpha}}}+\alpha_2^2 e^{-\frac{\mu t}{2}}\right)\cdot\left( \delta M_0^{\frac{3}{2}}M_1^2+M_0^{-1}M_2^2 \right) ||U||_{H^3}
				\\
				& \quad \quad+ C|\alpha_1\alpha_2| M_0^{-\frac{1}{2}}M_2^2 e^{-\frac{(2^{2-\alpha}\nu +\mu) t}{4}} ||\nabla B||_{H^3}.
			\end{split}
		\end{equation}
		\textbf{\\Step 3: Energy estimates of $U,B$.\\}
		Gathering above estimates in Step 1 and Step 2, we obtain
		\begin{equation*}
			\frac{1}{2}\frac{d}{dt}\left(||U||^2_{H^3}+||B||^2_{H^3}\right) + \frac{\nu}{2}||\nabla U||^2_{H^3}+ \frac{\mu}{2}||\nabla B||^2_{H^3}
			+ P(t)
			\leq  C \sum_{i=1}^{8}J_i,
		\end{equation*}
		where
		\begin{equation*}
			\begin{split}
			&P(t)=\kappa \sum_{|k|\leq 3} \int_{\mathbb{R}^3} \left( |(\nabla \times \Lambda^\alpha B) \times B|-|(\nabla \times \Lambda^\alpha B) \times \tilde{g}| \right)^2dx,
			\\
				&J_1=\left(||U||_{H^3}+||B||_{H^3}\right) \left( ||\nabla U||^2_{H^3}+ ||\nabla B||^2_{H^3}\right),
				\\
				& J_2=\left(||\tilde{f}||_{W{4,\infty}}+||\tilde{g}||_{W{4,\infty}}\right) \left(||U||^2_{H^3}+||B||^2_{H^3}\right),
				\end{split}
			\end{equation*}
			\begin{equation*}
			\begin{split}
				& J_3=\left( ||\tilde{f}\times \tilde{g}||_{H^3}+ M_0^{-\frac{1}{2}}M_2 |\alpha_2| e^{-\frac{\mu t}{4}} +|\alpha_1 \alpha_2| M_0^{-\frac{1}{2}} M_2^2 e^{-\frac{(2^{2-\alpha}\nu +\mu) t}{4}} \right) ||\nabla B||_{H^{3}},
				\end{split}
				\end{equation*}
				\begin{equation*}
				\begin{split}
				&J_4=\left(|\sigma| || (\nabla \times \tilde{g}) \times \tilde{g} ||_{H^3} +\kappa ||( (\nabla \times \tilde{g}) \times \tilde{g}) \times \tilde{g}  ||_{H^3} \right) ||\nabla B||_{H^{3}}
				\\
				&J_5=|\alpha_1|(M_1+M_2)e^{-\frac{\mu t}{4}} ||U||^2_{H^3}+ ||\tilde{g}||_{W^{5,\infty}}||B||^2_{H^3},
				\\
				& J_6=|\alpha_1| M_2 e^{-\frac{\mu t}{4}} ||B||^2_{H^3}+
				|\alpha_2|M_1 e^{-\frac{\nu t}{2^\alpha}}||U||_{H^3}||B||_{H^3},
				\\
				& J_7=\left( \alpha_1^2 e^{-\frac{\nu t}{2^{\alpha-1}}}+\alpha_2^2 e^{-\frac{\mu t}{2}}\right)\left( \delta M_0^{\frac{3}{2}}M_1^2+M_0^{-1}M_2^2 \right) ||U||_{H^3},
				\\
				& J_{8}= |\alpha_1\alpha_2| M_0^{-\frac{1}{2}}M_2^2 e^{-\frac{(2^{2-\alpha}\nu +\mu) t}{4}} ||\nabla B||_{H^3},
				\\
				& J_{9}=\left( ||\tilde{f}\times (\nabla \times \tilde{f})||_{H^3}+||\tilde{g}\times \left(\nabla \times \tilde{g}\right)||_{H^3} \right)||U||_{H^3},
				\\
				& J_{10}=\eta ||\nabla B||_{H^3} ||B||_{H^3} ||\tilde{g}||_{W^{3,\infty}}+M_0^{-\frac{1}{2}}M_2 |\alpha_1| e^{-\frac{\nu t}{2^\alpha}}  ||\Lambda^{\frac{\alpha}{2}} U||_{H^{3}},
				\\
				& J_{11}=C || B||^4_{H^3} \left( ||\tilde{g}||^2_{W^{4,\infty}}+||\tilde{g}||^2_{W^{2,6}} \right)+ C|| B||^2_{H^3}\left( ||\tilde{g}||^4_{W^{4,\infty}}+||\tilde{g}||^4_{W^{2,6}} \right).
			\end{split}
		\end{equation*}
		Using Lemma \ref{tuilun2}, Young's inequality and $\delta^{-\frac{1}{2}} \geq M_0 \gg 1$, we derive that
		\begin{equation}\label{327}
			\begin{split}
				&\frac{d}{dt}\left(||U||^2_{H^3}+||B||^2_{H^3}\right) + \left(\frac{\nu}{2}-C||U||_{H^3}-C||B||_{H^3} \right)||\Lambda^{\frac{\alpha}{2}} U||^2_{H^3}
				\\
				& \quad +\left(\frac{\mu}{2}-C||U||_{H^3}-C||B||_{H^3} \right)||\nabla B||^2_{H^3}
				\\
				\leq  & C\left( e^{-\frac{\nu t}{2^\alpha}} +e^{-\frac{\mu t}{4}} \right)\left(||U||^2_{H^3}+||B||^2_{H^3}+||B||^4_{H^3}\right)+C\left( M_0^{-1}+\delta^2 M_0^3\right) \left( e^{-\frac{\nu t}{2^\alpha}} +e^{-\frac{\mu t}{4}} \right)
			\end{split}
		\end{equation}
		for some constant $C$ depending on $M_1,M_2,\mu,\nu,\eta, \alpha_1,\alpha_2.$

		For $t \in [0,\infty)$, we assume that
		\begin{equation*}
			||U(t)||^2_{H^3}+||B(t)||^2_{H^3} \leq \frac{\min\{ \mu, \nu\}}{4C}.
		\end{equation*}
		In case $t=0$, the above estimate holds.
		Applying differential inequality \eqref{327}, Gronwall's inequality and $\delta \leq M_0^{-2}$, we have
		\begin{equation}\label{328}
			||U(t)||_{H^3}+||B(t)||_{H^3} \leq  M_0^{-\frac{1}{2}}.
		\end{equation}
		In a result,
		\begin{equation}\label{329}
			||U(t)||_{H^3}+||B(t)||_{H^3} \leq  M_0^{-\frac{1}{2}}
		\end{equation}
		for all $t \in [0,\infty)$. Therefore, we complete the proof of Theorem \ref{thm}.
	\end{proof}
	
	

\section*{Acknowlegement}
The author is supported by Education Department of Hunan Province, general Program(grant No. 17C0039); the State Scholarship Fund of China Scholarship Council (No. 201808430121) and Hunan Provincial Key Laboratory of Intelligent Processing of Big Data on Transportation, Changsha University of Science and Technology, Changsha; 410114, China.

	\medskip
	\medskip
	
\end{document}